\documentclass[11pt,reqno,draft]{amsart}
\usepackage{amssymb,verbatim}
 \let\Bbb\mathbb \let\Cal\mathcal
\let\le\leqslant \let\ge\geqslant \let\phi\varphi
 \let\TIL\widetilde \let\OVER\overline
\def\opn{\operatorname } 
\def\CC {{ \Bbb C {\,} }} \def\RR {{ \Bbb R {\,} }} \def\ZZ {{ \Bbb Z {\,} }}
\theoremstyle{definition}
\newtheorem*{DEF*}{Definition}

\theoremstyle{remark}
\newtheorem*{REM*}{Remark} 

\theoremstyle{plain}  
\newtheorem*{COR*}{Corollary}
\newtheorem*{LEM*}{Lemma}

\theoremstyle{plain}  
\newtheorem{PROP}{Proposition}[section]

\newtheorem{THM}{Theorem}[section]  

\numberwithin{equation}{section}

\newcommand{\diag}{\operatorname{diag}}

\def\vac{{1\!\!\Bbb I}}

\address{A.~M.~Vershik, St.~Petersburg Department of 
Steklov Institute of Mathematics, 27 Fontanka,
191023 St.~Petersburg, Russia.}
\email{vershik@pdmi.ras.ru}
\address{M.~I.~Graev, Institute for System Studies, 36-1
Nakhimovsky pr.,
117218 Moscow, Russia}
\email{graev\_36@mtu-net.ru}
\title[\today]{}

\begin{document}
\maketitle
\vskip-1cm
\vskip-1cm
\vskip-1cm


UDC 517.5


\begin{center}{}\sc
The basic representation of the current group
$O(n,1)^X$ in the $L^2$ space
\\
over the generalized Lebesgue measure
\\  \rm
A.~M.~Vershik\footnote{Supported by the grants 
NSh.-2251.2003.1 and RUMI-2662-ST-04.} and
M.~I.~Graev\footnote{Supported by RFBR, project 01-04-00363.}
\end{center}

\maketitle
\setcounter{tocdepth}{1}
{\footnotesize\tableofcontents}

\section{Introduction}

The current group $G^X$, where $X$ is a smooth manifold with a finite
continuous measure $m$ and $G$ is a Lie group, is the collection of
bounded piecewise continuous or, more generally, Borel 
$G$-valued functions on $X$. The group operation is 
the pointwise multiplication.

This paper continues the series of papers 
\cite{V-G-3, V-G-2, V-G-4, V-G-5, V-G-6, V-8}, in which irreducible unitary
representations of the current groups 
$G^X$ for some simple Lie groups $G$ were constructed that are invariant
under any $m$-preserving transformations of the space $X$: 
in the first paper
\cite{V-G-3} of this series, this was done
for $G=SL(2,\RR)$, and in
\cite{V-G-2,V-G-6}, for $G=SO(n,1)$ and $SU(n,1)$.
We called them the basic representations of the current groups
$G^X$ and described several models of these representations. Note
that basic representations are uniquely determined by their spherical
functions, which in
\cite{V-G-3} were called canonical states; see \cite{V-G-6}.

As shown in \cite{V-G-2} (see also \cite{V-9,V-10}), the series
$G=SO(n,1)$ and $SU(n,1)$ exhaust all simple Lie groups
$G$ for which such representations of the current groups do exist 
(in other words, groups with nontrivial first homology with
coefficients in some irreducible unitary representation).
Canonical states for a close group, the group of automorphisms
of an infinite homogeneous tree, were considered in
\cite{V-11}. See also \cite{O-15}, where canonical states on the group
$U(\infty )$ are considered.

In this paper, we present a detailed study of the so-called commutative model
of the basic representation of the group
$G^X$ for the series $G=O(n,1)$. We will return to the analysis of the series
$U(n,1)$ elsewhere. For convenience, we replaced the groups
$SO(n,1)$ with $O(n,1)$.

Among different realizations of the same unitary representation 
$\pi$ of an arbitrary metrizable topological group $G$, of special interest
are commutative models. Namely, let
$G_0 \subset G$ be a commutative subgroup such that the restriction of
$\pi$ to $G_0$ is cyclic. Then, by an isometric operator, one can
bring the operators $\pi(a)$,  $a\in G_0$, to diagonal form, i.e.,
realize the Hilbert space of the representation
$\pi$ as the space of square-integrable measurable functions with respect to the spectral
measure of the representation
$\pi|_{G_0}$, in which the representation operators corresponding to elements
of the subgroup $G_0$ act as multiplicators. We call such a realization
the commutative model of the representation $\pi$ with respect to the subgroup
$G_0$. The commutative models of irreducible unitary 
representations of the group
$SL(2,\RR)$ with respect to its commutative subgroups --- orthogonal,
unipotent, and diagonal --- are well known; see, for example, 
\cite{G-G-1}. Clearly, it makes sense to choose the subgroup
$G_0$ so that the representation operators corresponding to other elements of
the group $G$ will have simple expressions. 

In what follows, such subgroups are the unipotent subgroups
$Z$ and $Z^X$ of the groups
$O(n,1)$ and $O(n,1)^X$, respectively. The subgroup
$Z$ is commutative (in contrast to the maximal unipotent subgroup in
$U(n,1)$) and isomorphic to the additive group
$\RR^{n-1}$. The case
$n=2$, more exactly, that of the group
$SL(2,\RR)$, was earlier considered in
\cite{V-G-4,V-8}.

The interest to this realization is also
due to the fact that the diagonalization of operators
corresponding to the chosen commutative subgroup leads to remarkable
infinite $\sigma$-finite measures
$\nu$ in the space of distributions, which have a large group of linear
symmetries. In the case of
$SL(2,\RR)$, this measure $\nu$ turned out to be related to the 
well-known gamma process; see
\cite{V-G-4} and, for more details, \cite{V-8}. More exactly, this is a
$\sigma$-finite measure equivalent to the law of the gamma process and
invariant under the multiplication of realizations of the process
by functions with zero integral of the absolute value of logarithm 
(an analog of the linear transformation corresponding to 
a diagonal matrix with determinant $1$), which allows one
to call it the Lebesgue measure in the infinite-dimensional space.

In the general case
$G=O(n,1)$ considered in this paper, this measure $\nu$
is concentrated on vector distributions, and in addition to 
the above property it is also invariant under the pointwise action
of the group $O(n{-}1)^X$. In all cases, the measure
$\nu$ is concentrated on distributions that are linear combinations
of $\delta$-functions.

From the viewpoint of representation theory, the advantage of the
model constructed in
this paper is that the representation
operators have an explicit description, in contrast to the realization
in the Fock space
(see \cite{V-G-6}). Moreover, this model implicitly suggests certain 
advantages of the similar realization of unitary representations of
the Lie group $O(n,1)$ itself.

In the appendices we briefly discuss additional properties of the measures
arising in the construction of commutative models of representations and
describe the limiting case of the group
$O(\infty,1)$.

\section{Starting definitions and formulas}

\subsection{Current groups}

In this and the next sections we recall some definitions from
papers of the 60s and 70s
(see \cite{V-G-3,V-9,V-10,V-11,V-19,Guichardet} and references therein).

The current group $G^X$, where $X$ is a smooth manifold with a finite
continuous measure $m$ and $G$  is a Lie group, is the collection of 
bounded Borel $G$-valued functions on $X$. The group operation is the
pointwise multiplication. At first we do not introduce any topology on
$G^X$.

The problem is to construct irreducible unitary representations of current 
groups that are
invariant under $m$-preserving
transformations of the space $X$.
Such (invariant) representations are nonlocal, and one usually constructs
them in some or other realization of the Fock space, more
exactly, in a Hilbert space with factorization structure. 
It can be
not only the ordinary Fock space realized as the 
$L^2$ space over the standard Gaussian measure, but also the 
$L^2$ space over another measure corresponding to a L\'evy process,
i.e., a process constructed from an infinitely divisible distribution
on $\RR$ or $\RR^n$; for information on factorizations, see
\cite{V-19}. In this paper, such measures will be multidimensional analogs
of the classical gamma process and $\sigma$-finite measures
constructed from them.

\subsection{Canonical states}

A state on a topological group is a 
{\it positive definite continuous complex-valued Hermitian
($f(x^{-1})=\bar f (x)$) function normalized to the unity at the
group identity}; a state is called
{\it infinitely divisible} if it can be included into a continuous
one-parameter semigroup of states. In the representation-theoretic
language, the condition of infinite divisibility means
that the group representation corresponding to this state
(according to the GNS construction) can be included into an arbitrarily
high tensor power of some representation of the same group.

It is easy to show that 
the generator of the semigroup, i.e., the derivative of the one-parameter
family at the group identity, is a {\it conditionally positive definite
function on the group} (sometimes it is called a function of
negative type); at the same time, it is the squared norm of some
$1$-cocycle of the group with values in some 
unitary representation. This cocycle is cohomological to zero 
if and only if the generator, regarded as a function on the group,
is bounded in absolute value. In this case, it coincides,
up to sign and constant summand, with a positive
definite function.
That is why we are interested in unbounded
conditionally positive definite functions 
and the corresponding nontrivial $1$-cocycles with values in irreducible
unitary representations. Such cocycles do not always exist
and must lie in representations ``glued'' to the identity
representation.

{\it A canonical state on a group is an infinitely divisible state
such that the generator of the corresponding semigroup is
an unbounded conditionally positive definite function}.

Given a canonical state on a group $G$, we construct 
an irreducible unitary nonlocal representation (integral) of the corresponding
current group --- the group of $G$-valued measurable functions.
Conversely, each invariant (see above) representation
of the current group realized in the Fock space is generated by some
canonical state on $G$.

Thus a canonical state is the exponential of the squared norm 
of a nontrivial cocycle of the group with values in an irreducible
unitary representation of the group. The complete list of groups
for which a canonical state does exist is not yet known; however, there
is a number of examples, the most important of which are some
classical semisimple Lie groups, nontrivial $1$-cocycles on which were first found in
\cite{V-G-2,V-G-3}. Another example --- the group of automorphisms 
of a tree --- is considered in \cite{V-11}.

In the case of semisimple groups, one requires additionally that
a canonical state should be constant on a chosen maximal compact 
subgroup. This condition leads to a unique (up to the choice of
a positive degree) canonical state; namely, the corresponding 
one-parameter subgroup is the subgroup of spherical functions,
and its generator is the 
derivative of the family of spherical functions
of complementary series with respect to the parameter 
at the group identity.
Thus, among semisimple groups, canonical states and the corresponding 
nonlocal representations of current groups exist only on groups
of rank $1$ (and not all at that), namely, on 
$SO(n,1)$ and  $SU(n,1)$\footnote{It is possible 
to define a canonical state on groups
with the Kazhdan property (= the identity representation is 
an isolated point in the 
space of all irreducible unitary representations, as in a semisimple 
group of rank more than one), but in this case we must omit
the positivity condition and, consequently, allow considering
nonunitary representations. For example, for $SL(n, R)$,
$n>2$, the derivative of the family of spherical functions of 
complementary series representations makes sense, but
it will not be a conditionally positive function on the group.
Thus the exponential of this function  will not be a state on the group
in the true sense, but nevertheless it still determines an
invariant bilinear form on the group algebra, which is interesting for 
further construction of nonunitary representations of current 
groups. Note that F.~A.~Berezin \cite{Ber} considered
presumably similar notions,  but the relation between his construction 
and our notions is still not clear.}.

It is clear from above that the restriction of a canonical state to any
commutative subgroup is an infinitely divisible characteristic
(=positive definite) function, i.e., the Bochner--Fourier transform 
of an infinitely divisible {\it measure} on the group of characters of
the commutative subgroup; hence this measure determines a L\'evy 
process with values in the group of characters. In our case of the
group $SO(n,1)$, the commutative subgroup and its group of characters is
$\RR^{n-1}$, and the representation of the current group is realized in
the $L^2$ space over some vector L\'evy process.

Let $\phi(\cdot)$ be a canonical state on some group
$G$; considering its restriction to a commutative subgroup, 
we can write the Fourier transform of the law of the corresponding
L\'evy process as
$$\Phi(g(\cdot))=\exp\Bigl(\int_X -\log {\phi(g(x))}dx\Bigr).$$

\subsection{The group $G=O(n,1)$}
In this paper the group  $G=O(n,1)$ is realized as the collection of all
real matrices of order $n+1$ preserving the quadratic form
$$
2x_1x_{n+1}+x_2^2+ \ldots +x_n^2.
$$
In another formulation,
$G$ is the group of all real matrices
satisfying the relation
\begin{equation}{}\label{1-1}
gs g^*=s,\quad \text{где}\quad
s=\begin{pmatrix}0&0&1\\0&e&0\\1&0&0\end{pmatrix},
\end{equation}
where $e$ is the identity matrix of order
$n-1$ and $*$ stands for transposition.

We write elements $g\in G$ as block matrices
$$
g=\begin{pmatrix}{}g_{11}&g_{12}&g_{13}\\
g_{21}&g_{22}&g_{23}\\
g_{31}&g_{32}&g_{33}\end{pmatrix},
$$
where the diagonal blocks are quadratic matrices of orders
$1$, $n-1$, and $1$, respectively.

Condition \eqref{1-1} is equivalent to the following relations 
between the elements of these matrices:
\begin{gather*}{}
g_{13}g_{31}+g_{11}g_{33}+g_{12}g_{32}^*=1,\quad
g_{23}g_{21}^*+g_{21}g_{23}^*+g_{22}g_{22}^*=e,
\\
2g_{11}g_{13}+g_{12}g_{12}^*=0,\quad
2g_{31}g_{33}+g_{32}g_{32}^*=0,
\\
g_{11}g_{23}+g_{13}g_{21}+g_{22}g_{12}^*=0,\quad
g_{31}g_{23}+g_{33}g_{21}+g_{22}g_{32}^*=0.
\end{gather*}

The group $G$ contains as a subgroup the group
$Z$ of all block matrices of the form
$$
z=\begin{pmatrix}{}1&0&0\\
-\gamma^*&e&0\\-\frac{|\gamma|^2}{2}&\gamma&1\end{pmatrix}, \quad
\gamma \in\RR^{n-1},
$$
where $|\gamma | = (\sum {\gamma _i}^2)^{1/2}$.
The group $Z$ is commutative and isomorphic to the additive group 
$\RR^{n-1}$, and it is the maximal unipotent subgroup of
$G$. We identify elements
$z\in Z$ with vectors
$\gamma \in\RR^{n-1}$ and write $z(\gamma )$ or simply $\gamma $ instead of $z$.

By $D$ we denote the subgroup of block diagonal matrices 
from $G$,
and by $B$ the subgroup of all block lower triangular matrices.
Elements of the subgroup $D$ are matrices of the form
$$
d=\diag (\epsilon ^{-1},u, \epsilon ),\quad \epsilon \in\RR\setminus 0,\quad
u\in O(n{-}1).
$$
The group $B$ is the semidirect product
$B=Z\leftthreetimes D$, and the group
$D$ is the direct product of the groups $\RR^*$ and $O(n{-}1)$.

The canonical state $\phi(g)$ on the group $G$ is given
by the following formula:
$$
\phi(g)=\exp\Bigl(\Bigl.\frac{d\psi_{\lambda }(g)}{d \lambda }\Bigr|_{\lambda =0}\Bigr),
$$
where $\psi_{\lambda }$ is the spherical function of the 
complementary series representation of the group $G$
with parameter $\lambda$, the value
$\lambda =0$ corresponding to the special representation of $G$ glued to
the identity representation; see
\cite{V-G-3} and \cite{V-G-6}.

Note that the special representation of the group $G$ corresponding to
$\lambda =0$ has a nontrivial $1$-cocycle and is the unique irreducible
unitary representation of $G$ with this property.

\subsection{The current group $G^X=O(n,1)^X$}

We consider the group $G=O(n,1)$ in the realization described above.
Denote by $Z^X$, $D^X$, and $B^X$ the subgroups of functions $g(x)\in G^X$
with values in the subgroups $Z$, $D$, and $B$, respectively. Note that
the group $Z^X$ is isomorphic to the additive group
$(\RR^{n-1})^X$ of bounded Borel functions on
$X$ with values in $\RR^{n-1}$, and the group $B^X$ 
is the semidirect product
$B^X=Z^X\leftthreetimes D^X$.

In \cite{V-G-2,V-G-6}, for each group $O(n,1)^X$,
an irreducible unitary representation was constructed that is invariant
under $m$-preserving transformations of the space $X$.
It was called the basic representation; several models
of this representation were described.

In the chosen realization of the
group $O(n,1)$,
the spherical function $\Psi(g(\cdot))$ (canonical state) on the current group
$G^X$ that determines this representation
is given by the following formula:
\begin{equation}{}\label{1-11}
\Psi(g(\cdot))=\exp\Bigl(-\frac12\int_X\log
\Big|
\frac{g_{11}(x)+g_{33}(x)-
g_{13}(x)-g_{31}(x)}{2}
\Big|
\,d\,m(x)\Bigr),
\end{equation}
where $g_{ij}$ are elements of the block matrix $g\in G^X$.
Its restriction to the infinite-dimensional commutative current group
$Z^X$ is given by the formula
\begin{equation}{}\label{1-12}
\Psi(\gamma (\cdot))=\exp\Bigl(-\frac12\int_X
\log\Bigl(1+\frac{|\gamma(x)|^2}{4}\Bigr)
\,d\,m(x)\Bigr),
\end{equation}
where $\gamma $ is the function determining the block matrix 
$$
z=\begin{pmatrix}{}1&0&0\\
-\gamma^*&e&0\\-\frac{|\gamma|^2}{2}&\gamma&1\end{pmatrix}\in Z^X.
$$
Automatically,
this restriction is the characteristic function of an
infinitely divisible distribution, because the canonical state on the group $G$
is infinitely divisible.

We construct a commutative model of the basic representation of the group
$O(n,1)^X$. It is realized in the Hilbert space
$L^2(\nu)$ of functions on the space of vector distributions
square-integrable with respect to the measure $\nu$ introduced below.
By the properties of this measure, on
$L^2(\nu)$ there is a natural unitary representation 
$U_g$ of the block triangular subgroup
$B^X$, where elements $z\in Z^X$ act as diagonal operators.
The desired commutative model of the basic representation of the group
$G^X$ is obtained by extending this representation from the subgroup
$B^X$ to the whole group $O(n,1)^X$.

The measure $\nu$ is determined by its density with respect to another measure
$\mu $ introduced in this paper. By definition,
$\mu$ is the finite measure on the space of vector distributions whose
Fourier transform is the functional
$\Psi(\gamma (\cdot))$ given by \eqref{1-12}.

Note that the measure $\mu $ is the law of the L\'evy process
obtained by the canonical construction from the infinitely divisible measure 
$\alpha$ on the additive group $\RR^{n-1}$ whose Fourier transform is
$(1+\frac{|\gamma |^2}{4})^{1-n}$,
$\gamma \in\RR^{n-1}$. The density of this measure $\alpha$
with respect to the standard Lebesgue measure is given, up to constant factor,
by the following formula:
$$
\TIL\psi(\xi )=|\xi |^{\frac{n-1}{2}}\,K_{\frac{n-1}{2}}(2|\xi|),
$$
where $K_{\nu}(x)$ is the modified Bessel function of the third kind, see
\cite{B-12}. In particular, 
$\TIL\psi(\xi )   =   e^{-2|\xi|}$ in the case $n=2$. The measure
$\alpha $ is $O(n{-}1)$-invariant, and it is a multidimensional generalization of
the symmetrized gamma distribution.

\subsection{Several integral relations} 
In further constructions we will use formulas for the Fourier transform
of the functions
$(1+\frac{|\gamma|^2}{4})^{-\lambda/2}$ and 
$|\gamma|^{-\lambda}$ on $\RR^{n-1}$, where
$|\gamma|=(\sum_{i=1}^{n-1}\gamma_i^2)^{1/2}$. Namely, 
\begin{equation}{}\label{17-3}
\int_{\RR^{n-1}} \Bigl(1+\,\frac{|\gamma|^2}{4}\Bigr)^{-\lambda/2}\,e^{\,i\,<\xi,\gamma>}\,d\gamma =
c_n\frac{2}{\Gamma (\lambda /2)}\,|\xi|^{(\lambda -n+1)/2}\,
K_{(n-1- \lambda )/2}(2|\xi|),
\end{equation}
\begin{equation}{}\label{17-4}
\int_{\RR^{n-1}} |\gamma|^{-\lambda}\,e^{\,i\,<\xi,\gamma>}\,d\gamma =
c_n\frac{2^{\,-\lambda}\, \Gamma ((n-1-\lambda)/2) }{\Gamma (\lambda /2)}
\,|\xi|^{\lambda -n+1},
\end{equation}
where the coefficient $c_n$, the same in
\eqref{17-3} and \eqref{17-4}, depends only on $n$.

The Bessel function $K_{\rho}$ is given by the 
following equation (see \cite{B-12}):
\begin{equation}{}\label{17-5}
K_{\rho}(2z)=\frac{\pi}{2\sin(\pi\rho)}\,(I_{-\rho}(2z)-I_{\rho}(2z)),
\end{equation}
where
\begin{equation}{}\label{17-6}
I_{\rho}(2z)=\sum_{m=0}^\infty \frac{z^{2m+\rho}}{m! \,\Gamma (m+\rho+1)}.
\end{equation}

It is important for the sequel that this function
is continuous and strictly positive on the half-line
$0<x< \infty$.

\begin{REM*}{}
For integer values of $\rho$, the series for $K_{\rho}$ contains terms with
$\log z$; for half-integer values of $\rho$, the expression for $K_{\rho}$ 
can be simplified:
$$
K_{n+\frac12}(z)=\Bigl(\frac{\pi}{2z}\Bigr)^{1/2}\,e^{-z}\,
\sum_{k=0}^n \frac{(n+k)!}{k!(n-k)!(2z)^k}.
$$
\end{REM*}

It is convenient to write equation \eqref{17-3} in another form. Namely,
setting $\lambda =2 \rho +n-1$, we obtain
\begin{equation}{}\label{717-3}
\int_{\RR^{n-1}} \Bigl(1+\frac{|\gamma|^2}{4}\Bigr)^{-\frac{n-1}{2}-\rho}\,e^{\,i\,<\xi,\gamma>}\,d\gamma =
c_n\frac{2}{\Gamma (\frac{n-1}{2}+\rho)}\,|\xi|^{\rho}\,
K_{\rho}(2|\xi|).
\end{equation}

Let us give a brief derivation of formulas
\eqref{17-3} and \eqref{17-4} for $n>2$.
In spherical coordinates
the first integral takes the form
$$
J_1=c\,\int_0^{\infty }\int_0^{\phi} \Bigl(1+\frac{r^2}{4}\Bigr)^{- \lambda /2}\,
e^{\,i|\xi|\cos\phi}\,r^{n-2}\,\sin^{n-3}\phi\,d\phi\,dr.
$$
Integrating with respect to $\phi$ yields (see \cite[formula 3.915.5]{G-R-15})
$$
J_1=c\,|\xi|^{-\frac{n-3}{2}}\int_0^{\infty }
\Bigl(1+\frac{r^2}{4}\Bigr)^{- \lambda /2}\,r^{\frac{n-1}{2}}\,
J_{\frac{n-3}{2}}(|\xi|r)\,dr,
$$
where $J_{\rho}(\gamma)$ is the Bessel function of the first kind
(see \cite{B-12}). Similarly, for the integral
\eqref{17-4} we obtain
$$
J_2=c\,|\xi|^{-\frac{n-3}{2}}\int_0^{\infty }
r^{\frac{n-1}{2}- \lambda }\,
J_{\frac{n-3}{2}}(|\xi|r)\,dr.
$$

Integrating with respect to $r$ yields the expressions
\eqref{17-3} and \eqref{17-4} for
$J_1$ and $J_2$, respectively; see \cite[formulas 6.565.4 and 6.561.14]{G-R-15}.

\begin{THM}[Multidimensional analog of the L\'evy--Khintchin formula]{}\label{THM:345-1}
The function $\log(1+\frac{|\gamma|^2}{4})$ on $\RR^{n-1}$
has the following integral representation:
\begin{equation}{}\label{345-1}
\log\Bigl(1+\frac{|\gamma |^2}{4}\Bigr)=
\int_{\RR^{n-1}}\Bigl(e^{\,i<\xi, \gamma >} - 1\Bigr)\,
g(\xi)\,d\xi, \quad \text{where}\quad g(\xi)=|\xi|^{-\frac{n-1}{2}}\,
K_{-\frac{n-1}{2}}(2|\xi|).
\end{equation}
\end{THM}
Indeed, it follows from \eqref{17-3} that
$$
\int_{\RR^{n-1}}
\Bigl(e^{\,i<\xi, \gamma >}-1\Bigr)\,|\xi|^{\frac{\lambda-n+1}{2} }\,
K_{\frac{\lambda-n+1}{2}}(2|\xi|))^{1/2}\,d\xi
= c\, \Gamma (\lambda /2)\,
\Bigl(\Bigl(1+\frac{|\gamma |^2}{4}\Bigr)^{-\lambda /2}-1\Bigr).
$$
As $\lambda \to 0$, we obtain \eqref{345-1}.

\section{The commutative model of the complementary series of
irreducible unitary representations of the group
$O(n,1)$}

\subsection{The action of the group $G$ on $\RR^{n-1}$ and the
$1$-cocycle}
Let  $Y$ be the manifold of one-dimensional subspaces in
$\RR^{n+1}$ lying inside the light cone
$$
2x_1x_{n+1}+x_2^2+ \ldots +x_n^2=0.
$$
The group $G=O(n,1)$, regarded as a group of linear transformations in 
$\RR^{n+1}$, acts transitively on
$Y$. We use the right notation for this action:
$y\to y\bar g$.  Note that in another interpretation 
$Y$ is the absolute of the $n$-dimensional Lobachevsky space 
realized as the collection of one-dimensional subspaces in 
$\RR^{n+1}$ lying inside the light cone.

Let us realize $Y\setminus y_0$, where $y_0=(\lambda ,0, \dots ,0)$,
as the intersection of the cone with the hyperplane 
$x_{n+1}=1$, i.e., as the set of points in 
$\RR^{n+1}$ of the form
$$
\Bigl(-\frac{|\gamma|^2}{2}, \gamma _i, \dots , \gamma _{n-1},1\Bigr),
$$
where $\gamma =(\gamma _i, \dots , \gamma _{n-1})\in\RR^{n-1}$ and
$|\gamma | = (\sum {\gamma _i}^2)^{1/2}$.  According to this realization, there
is a natural bijection $Y\setminus y_0\to\RR^{n-1}$; hence the action of the 
group $G$ on $Y$ induces an action $\gamma \to \gamma \bar g$ of $G$
on the space $\RR^{n-1}$. We emphasize that this action is not linear.

It follows from the definition that the vector
$\gamma \bar g$ is given by the following formula:
\begin{equation}{}\label{1-2}
\gamma \bar g=\Bigl(-\frac{|\gamma|^2}{2}g_{13}+ \gamma g_{23}+g_{33}\Bigr)^{-1}\,
\Bigl(-\frac{|\gamma|^2}{2}g_{12}+ \gamma g_{22}+g_{32}\Bigr),
\end{equation}
where $g_{ij}$ are elements of the block matrix $g$. In particular,
$$
\begin{array}{ll}
\gamma\OVER g= \gamma + \gamma_0   &\quad \text{for}\quad g=z(\gamma _0)\in Z;
\\
\gamma\OVER g= \epsilon^{-1} \gamma u &\quad \text{for}\quad
g=\diag(\epsilon ^{-1},u, \epsilon );
\\
\gamma\OVER s= -\dfrac{2 \gamma }{|\gamma |^2}.&
\end{array}
$$

Now let us define a function $\beta (\gamma,g)$ by the formula
\begin{equation}{}\label{1-3}
\beta (\gamma ,g)=\Bigl|-\frac{|\gamma|^2}{2}g_{13}+ \gamma g_{23}+g_{33}\Bigr|, \quad
\gamma \in\RR^{n-1}, \quad g\in G.
\end{equation}
In particular,
$$
\begin{array}{ll}
\beta (\gamma ,g)=1&\quad \text{for}\quad g\in Z;
\\
\beta (\gamma ,g)= |\epsilon\,|& \quad\text{for}\quad g=\opn{diag}(\epsilon ^{-1},u, \epsilon );
\\
\beta (\gamma ,s)=\dfrac{|\gamma |^2}{2}.&
\end{array}
$$
It follows from the definition that $\beta (\gamma ,g)$ is a $1$-cocycle
of the group $G$ with values in $\RR^*$, i.e.,
\begin{equation}{}\label{1-4}
\beta (\gamma ,g_1g_2)=\beta (\gamma ,g_1)\,\beta (\gamma\bar g_1,g_2) \quad
\text{for any} \quad \gamma \in\RR^{n-1}\quad \text{and}\quad  g_1,g_2\in G.
\end{equation}

\subsection{The standard model of the complementary series representations}

Each irreducible unitary complementary series representation 
of the group $G=O(n,1)$ is determined by a number $\lambda $ from 
the interval $0< \lambda <n-1$. In the standard realization, the representation
$T^{\lambda }$ with parameter $\lambda$ acts in the Hilbert space
$L_{\lambda }$ of real-valued functions $f(\gamma )$ on
$\RR^{n-1}\simeq Z$ with scalar product
\begin{equation}{}\label{1-70}
<f_1,f_2>\,\,=\!\!\!\!\!\int\limits_{\RR^{n-1} \times \RR^{n-1}}\!\!\!\!\!\!\!%
| \gamma '- \gamma ''|^{- \lambda }
\,f_1(\gamma ')\,f_2(\gamma '')\,d  \gamma '\,d \gamma '',
\end{equation}
where $d \gamma =d \gamma _1 \ldots d \gamma _{n-1}$ is the Lebesgue measure on
$\RR^{n-1}$. The operators of this representation have the form
\begin{equation}{}\label{1-8}
T^{\lambda }_g f(\gamma )=f(\gamma \bar g)\,
\beta ^{1-n+\frac{ \lambda}{2} } (\gamma ,g),
\end{equation}
where $\gamma \bar g$ and $\beta (\gamma ,g)$ are given by equations
\eqref{1-2} and \eqref{1-3}, respectively. In particular,
\begin{align}{}
\label{273}
&
T^{\lambda }_{z}f(\gamma )=f(\gamma + \gamma _0)\quad \text{for}\quad
z = z(\gamma _0) \in Z;
\\ &
\label{274}
T^{\lambda }_d f(\gamma )= |\epsilon|^{1-n+ \frac{\lambda}{2} }\,
f(\epsilon ^{-1} \gamma u)\quad \text{for}\quad
d=\diag(\epsilon  ^{-1},u, \epsilon );
\\ &
\label{275}
T^{\lambda }_s f(\gamma ) =  f\Bigl(-\frac{2 \gamma }{|\gamma|^2}\Bigr)\,
\Bigl( \frac{|\gamma|^{2}}{2} \Bigr)^{1-n+ \frac{\lambda}{2} }
\quad \text{for}\quad
s=\begin{pmatrix}{}\,0&0&1\,\\\,0&e&0\,\\\,1&0&0\,\end{pmatrix}.
\end{align}

The group property 
of these operators follows immediately from the property
\eqref{1-4} of the function $\beta (\gamma ,g)$,
and their unitarity follows from the relations
\begin{equation}{}\label{1-5}
d(\gamma \bar g)=\beta^{1-n} (\gamma,g)\,d \gamma \quad
\text{for any} \quad g\in G,
\end{equation}
where $d \gamma =d \gamma _1 \ldots d \gamma _{n-1}$, and
\begin{equation}{}\label{1-6}
|x-y|^2=|x\bar g-y\bar g|^2\, \beta (x,g)\,\beta (y,g)
\end{equation}
for any $x,y\in\RR^{n-1}$ and $g\in G$.

Relations \eqref{1-5} and \eqref{1-6} are easily verified for elements from
$Z$ and $D$ and the element $s$. It follows from the properties of the 
$1$-cocycle $\beta (\gamma ,g)$ that they hold for any element
$g\in G$.

\subsection{Construction of the commutative model of complementary
series representations}

Let us describe the commutative model of a complementary series
representation $T^{\lambda }$ of the group
$G=O(n,1)$ with respect to the subgroup $Z$, i.e., the model in which the operators
$T^{\lambda }_z$,  $z\in Z$, act as multiplicators.

This model is obtained by passing from functions
$f(\gamma )$ in the standard model to their Fourier transforms
$$
\phi(\xi)=\int_{\RR^{n-1}} e^{\,i<\xi, \gamma >} f( \gamma )\,d \gamma.
$$

\begin{THM}{}\label{THM:1} In the commutative model, the 
complementary series representation
$T^{\lambda }$ is realized in the Hilbert space
$L_{\lambda }$ of complex-valued functions on
$\RR^{n-1}$ with the norm
\begin{equation}{}\label{1-13}
\|\phi\|^2=
\tfrac{2^{\,-\lambda}\, \Gamma ((n-1-\lambda)/2) }{\Gamma (\lambda /2)}\,
\int_{\RR^{n-1}} |\xi|^{1-n+ \lambda }\,|\phi(\xi)|^2\,d\xi,
\quad
|\xi|=\,<\xi,\xi>^{1/2},
\end{equation}
satisfying the condition
\begin{equation}{}\label{1-133}
\OVER{\phi(\xi)}=\phi(-\xi).
\end{equation}
The operators of the representation are given by the formula
\begin{equation}{}\label{1-14-1}
T^{\lambda }_g \phi(\xi)=\int_{\RR^{n-1}} A^{\lambda }(\xi,\xi',g)\,
\phi(\xi')\,d\xi',
\end{equation}
where
\begin{equation}{}\label{1-15-1}
A^{\lambda }(\xi,\xi',g)=\int_{\RR^{n-1}}
e^{\,i\,(<\xi,\gamma>-<\xi',\gamma\OVER g>)}\,
\beta ^{1-n+\lambda/2 }(\gamma,g)\,d\gamma.
\end{equation}
In particular, 
$$
\begin{array}{ll}T^{\lambda }_{z}\phi(\xi)=
e^{-i<\xi,\gamma_0>}\,\phi(\xi)&\quad \text{for}\quad
z = z(\gamma _0) \in Z;
\\
T^{\lambda }_d \phi(\xi)= |\epsilon|^{\lambda/2 }\,
\phi(\epsilon\, \xi u)&\quad \text{for}\quad
d=\diag(\epsilon  ^{-1},u, \epsilon ).
\end{array}
$$
\end{THM}

\begin{proof}{} In the new model, the squared norm is given by the formula
$$
\|\phi\|^2=\!\!\!\!\!\!\!\!%
\int\limits_{\RR^{n-1} \times \RR^{n-1}}\!\!\!\!\!\!\!\!%
R(\xi,\xi')\,\phi(\xi)\,\OVER{\phi(\xi')}\,d\xi\,d\xi',
$$
where
$$
R(\xi,\xi')=\!\!\!\!\!\!\!\!%
\int\limits_{\RR^{n-1} \times \RR^{n-1}}\!\!\!\!\!\!\!%
|\gamma - \gamma '|^{- \lambda }
e^{\,i(<\xi, \gamma >-<\xi', \gamma' >)}\,d \gamma \,d\xi'=
\delta (\xi-\xi')\!\!\!\int\limits_{\RR^{n-1}} |\gamma |^{- \lambda }\,
e^{\,i<\xi, \gamma >}\,d\xi.
$$
This implies \eqref{1-13} in view of \eqref{17-4}.

Relation \eqref{1-133} is equivalent to the condition that 
the original space is real.

The formulas for representation operators in the new model
can be obtained directly from the formulas for these operators in
the original model by passing from functions
$f(\gamma )$ to their Fourier transforms.
\end{proof}

\begin{PROP}{}\label{PROP:3} In the commutative model of the representation
$T^{\lambda }$, the kernel
$A(\xi,\xi')=A^{\lambda }(\xi,\xi',s)$
of the operator $T^{\lambda }_s$ corresponding to the element
$s=\begin{pmatrix}{}0&0&1\\0&e&0\\1&0&0\end{pmatrix}$ has the following form:
\begin{eqnarray}{}\label{73-1}
A(\xi,\xi')&=&2^{1-\frac{\lambda }{2}}\int_0^{\infty }
\cos\Bigl(\xi x+\frac{2\xi'}{x}\Bigr)\,x^{ \lambda-2 }\,dx \quad \text{for}\quad n=2,\\
\label{63-1}
A(\xi,\xi')&=&c_n 2^{- \lambda/2 }\,
\int_0^\infty r^{\lambda -n}\,|r\xi+2r^{-1}\xi'|^{-\frac{n-3}{2}}\,
J_{\frac{n-3}{2}}(|r\xi+2r^{-1}\xi|)\,dr
\end{eqnarray}
for $n>2$, where $J_{\frac{n-3}{2}}$ is the Bessel function of the first kind.
\end{PROP}

Indeed, since $\gamma \OVER s=\dfrac{-2\gamma \,\,}{|\gamma |^2}$ and
$\beta (\gamma ,s)=\dfrac{|\gamma|^2}{2}$, it follows from 
\eqref{1-15-1} that
$$
A(\xi,\xi')
=  2^{n-1- \lambda/2 }\,\int_{\RR^{n-1}} e^{\,i\,(<\xi, \gamma >
                +<\xi',\frac{2 \gamma }{ |\gamma|^2}>)
}\,
|\gamma|^{2-2n+\lambda }\,d \gamma .
$$
For $n=2$, \eqref{73-1} follows immediately. 
For $n>2$, in spherical coordinates we obtain
$$
A(\xi,\xi)=c_n 2^{- \lambda/2 }\,
\int_0^\infty \int_0^{\pi}
e^{\,i\,|r\xi+2r^{-1}\xi'|\cos\phi}\,r^{\lambda -n}\,
\sin^{n-3}\phi \,d\phi\,dr.
$$
Integrating with respect to $\phi$ yields \eqref{63-1}.

\begin{REM*}{} For $n=2$, the kernel $A(\xi,\xi')$ can be expressed
in terms of Bessel functions:
\begin{eqnarray*}{}
A(\xi,\xi')&=&c\,\Bigl(\cos\frac{\pi \lambda }{2}\Bigr)^{-1}\,
|\xi'\xi^{-1}|^{1/2}\,[J_{\lambda -1}(2^{3/2}|\xi\xi'|)-
J_{1-\lambda }(2^{3/2}|\xi\xi'|)] \quad \text{for} \quad \xi\xi'<0,
\\
A(\xi,\xi')&=&c\,\Bigl(\cos\frac{\pi \lambda }{2}\Bigr)^{-1}\,
|\xi'\xi^{-1}|^{1/2}\,[I_{\lambda -1}(2^{3/2}|\xi\xi'|)-
I_{1-\lambda }(2^{3/2}|\xi\xi'|)] \quad \text{for} \quad \xi\xi'>0.
\end{eqnarray*}
\end{REM*}

The representation $T^{\lambda }$ 
is uniquely determined by its spherical function
$$
\psi_{\lambda }(g)=\,<T^{\lambda }_g\vac_{\lambda },\vac_{\lambda }>,
$$
where $\vac_{\lambda }\in L_{\lambda }$ is a vector of norm $1$ 
that is invariant
under the maximal compact subgroup of $G$ (vacuum vector).
In the chosen realization of the group
$G$, this spherical function takes the form
\begin{equation}{}\label{521-1}
\psi_{\lambda }(g) = \Bigl|\frac{g_{11}(x)+g_{33}(x)-
g_{13}(x)-g_{31}(x)}{2}\Bigr|^{-\lambda /2},
\end{equation}
where $g_{ij}$ are elements of the block matrix
$g$. In particular,
$$
\psi_{\lambda }(z(\gamma ))=\Bigl(1+\frac{|\gamma |^2}{4}\Bigr)^{-\lambda /2}.
$$

Let us introduce the following vector in the space
$L_{\lambda }$:
\begin{equation}{}\label{492}
f_{\lambda }(\xi) = \bigl(|\xi|^{-\frac{\lambda-n+1}{2} }\,
 K_{\frac{\lambda-n+1}{2}}(2|\xi|)\bigr)^{1/2},
\end{equation}
where $K_{\rho}(x)$ is the Bessel function defined above.

\begin{PROP}{}\label{PROP:432}  The vector $f_{\lambda }$ is invariant
under the maximal compact subgroup of the group $G$, hence it is
proportional to the vacuum vector $\vac_{\lambda }$.
\end{PROP}

\begin{COR*}{} The following equation holds:
$$
<T^{\lambda }_{z{\gamma }} f_{\lambda },f_{\lambda }>\,=
\|f_{\lambda }\|^2\,\psi_{\lambda }(z(\gamma )),
$$
i.e.,
\begin{equation}{}\label{342-1}
\int_{\RR^{n-1}}
e^{\,i<\xi, \gamma >}\,|\xi|^{\frac{\lambda-n+1}{2} }\,
K_{\frac{\lambda-n+1}{2}}(2|\xi|))^{1/2}\,d\xi=
\|f_{\lambda }\|^2\,\Bigl(1+\frac{|\gamma |^2}{4}\Bigr)^{-\lambda /2}.
\end{equation}
\end{COR*}

\begin{PROP}{}\label{PROP:342-2} We have $\|f_{\lambda }\|^2 =
c \Gamma (\lambda /2)\,\Bigl(1+\frac{|\gamma |^2}{4}\Bigr)^{-\lambda /2}.$
\end{PROP}

Indeed, it follows from \eqref{17-3} that the left-hand side of 
\eqref{342-1} is equal to $(2 c_n)^{-1}\, \Gamma (\lambda /2)$.

\subsection{The embedding $L_{\lambda }\to \bigotimes _{i=1}^l L_{\lambda _i}$}

\begin{PROP}{}\label{PROP:5} For any positive real numbers
$\lambda _1, \dots ,\lambda _l$,
$\sum \lambda _i<n-1$, there exists an isometric embedding
$$
\tau:  L_{\lambda }\to \bigotimes _{i=1}^l L_{\lambda _i},\quad
\lambda =\sum \lambda _i,
$$
that commutes with the action of the group
$G$. In the standard realization of representations it is given by the formula
\begin{equation}{}\label{31-21-1}
\tau f(\gamma^1, \dots ,\gamma^l)=f(\gamma^1)
\prod_{i=2}^l \delta (\gamma ^1- \gamma^i),
\end{equation}
and in the commutative model, by the formula
\begin{equation}{}\label{31-21-2}
\tau \phi(\xi _1, \dots ,\xi _l)=\phi(\xi _1+ \ldots +\xi_l).
\end{equation}
\end{PROP}

\begin{proof}{} First let us consider the standard model of representations.
Let $<\,,\,>$ and $<\,,\,>_l$ be the scalar products in the spaces
$L_{\lambda }$ and $\bigotimes _{i=1}^l L_{\lambda _i}$, respectively.
Then it is obvious that
$$
<\tau f,\tau f>_l\,\,=\!\!\!\!\!\!\int\limits_{\RR^{n-1} \times \RR^{n-1}}
\prod_{i=1}^l| \gamma '- \gamma ''|^{- \lambda_i }
\,f(\gamma ')\,f(\gamma '')\,d  \gamma '\,d \gamma ''=\,<f,f>.
$$
Thus the mapping $\tau$ is isometric. Obviously, it commutes with 
the action of $G$.

Now let us consider the commutative model of representations.
Let $\phi(\xi)$ be the Fourier transform of a function
$f(\gamma )$. Then, according to
\eqref{31-21-1}, the image of $\phi$ under this embedding is equal to
$$
\int f(\gamma _1)\,\prod_{i=2}^l
\delta  (\gamma _1- \gamma _i)\,\prod_{i=1}^l e^{\,i\,<\xi_i, \gamma_i> }\,d\gamma _i
=\phi(\xi _1+ \ldots +\xi_l).
$$
\end{proof}

\subsection{The commutative model of the special representation of the group
$O(n,1)$}

The special representation of the group
$O(n,1)$ is the irreducible unitary representation of this group glued to 
the identity representation. It is obtained from the complementary
series representations in the
$\lambda \to 0$ limit. Thus Theorem~\ref{THM:1} implies
the following result.

\begin{THM}{}\label{THM:111} In the commutative model, the special
representation 
$T^0$ of the group $O(n,1)$ is realized in the Hilbert space of complex-valued
functions on $\RR^{n-1}$ with the norm
\begin{equation}{}\label{11-13}
\|\phi\|^2 = \int_{\RR^{n-1}} |\xi|^{1-n}\,|\phi(\xi)|^2\,d\xi
\end{equation}
satisfying the condition $\OVER{\phi(\xi)}=\phi(-\xi)$.
Operators of the representation are given by the formula
\begin{equation}{}\label{2-14-1}
T_g \phi(\xi)=\int_{\RR^{n-1}} A^0(\xi,\xi',g)\,
\phi(\xi')\,d\xi',
\end{equation}
where
\begin{equation}{}\label{2-15-1}
A^0(\xi,\xi',g)=\int_{\RR^{n-1}}
e^{\,i\,(<\xi,\gamma>-<\xi',\gamma\OVER g>)}\,
\beta ^{1-n }(\gamma,g)\,d\gamma.
\end{equation}
In particular, 
\begin{eqnarray*}{}T_{z}\phi(\xi)&=&
e^{-i<\xi,\gamma_0>}\,\phi(\xi)\quad \text{for}\quad
z = z(\gamma _0) \in Z;
\\
T_d \phi(\xi)&=& \phi(\epsilon\, \xi u)\qquad\qquad \text{for}\quad
d=\diag(\epsilon  ^{-1},u, \epsilon ).
\end{eqnarray*}
\end{THM}

The special representation $T^0$ has a nontrivial $1$-cocycle
$\beta : G\to L_0$, where $L_0$ is the space of $T^0$.
Namely, set $f_0(\xi)=\lim_{\lambda \to 0} f_{\lambda }(\xi)$, where
$f_{\lambda }(\xi)\in L_{\lambda }$ is given by \eqref{492}, i.e.,
$$
f_0(\xi)=(|\xi|^{\frac{n-1}{2} }\, K_{\frac{n-1}{2}}(2|\xi|))^{1/2}.
$$
This vector is invariant under the maximal compact subgroup 
of $G$, but it does not belong to the space
$L_0$ of the special representation. The desired nontrivial 
$1$-cocycle  $G\to L_0$ is given by the equation
$$
\beta (\gamma ,g)=T^0_g f_0(\xi)-f_0(\xi ).
$$
It is known that $T^0$ is the unique irreducible unitary representation
of the group $G$ possessing a nontrivial $1$-cocycle.

\section{The measures $\mu$ and $\nu$ on the space of vector distributions}

\subsection{The spaces $ F$ and $\Phi$} 

The construction of the commutative model of the basic representation
of the current groups $G^X$ will be based on two remarkable
measures in the space of vector distributions.

Denote by $ F$ the space of smooth bounded real-valued functions on $X$
and by $\Phi$ the dual space with the ordinary topology (the space
of distributions on $X$). We will denote by
$<\,,\,>$ the pairing of elements from $ F$ and $\Phi.$

Let us introduce the space $ F^{n-1}$ of vector functions
$\gamma (x)=(\gamma_1 (x), \dots ,\gamma_{n-1} (x))$, $\gamma_i\in F$,
and denote by $\Phi^{n-1}$ the dual space. Elements from
$\Phi^{n-1}$ are realized as vector distributions
$\xi(x)=(\xi_1 (x), \dots ,\xi_{n-1} (x))$, $\xi_i\in\Phi$, with
the pairing
$$
<\xi,\gamma>\,=\sum_{k=1}^{n-1} <\xi_k,\gamma_k>.
$$
For example, if $\xi=c\, \delta _{x_0}$, where г$c=(c_1, \dots ,c_{n-1})
\in\RR^{n-1}$ and $\delta _{x_0}$ is a $\delta$-function on $X$, then
$<\xi,\gamma>\,=\sum_{k=1}^{n-1} c_k\, \gamma_k (x_0).$

Since in what follows we will consider only the spaces
$ F^{n-1}$ and
$\Phi^{n-1}$, we will omit the index $n-1$ in their notation.

\subsection{The measure $\mu$ on $\Phi$}
Let us introduce the following function on $\RR^{n-1}$:
\begin{equation}{}\label{181-1}
l(\gamma )=\Bigl(1+\frac{|\gamma |^2}{4}\Bigr)^{-1/2},\quad  \gamma \in\RR^{n-1}.
\end{equation}
It is known that this function is positive definite. According to
\eqref{17-3}, its Fourier transform equals
$$
\int_{\RR^{n-1}} \Bigl(1+\frac{| \gamma | ^2}{4}\Bigr)^{-1/2}\,e^{\,i\,<\xi,\gamma>}\,d\gamma =
c_n\frac{2}{\Gamma (1 /2)}\,|\xi|^{(-n+2)/2}\,
K_{(n-2)/2}(2|\xi|).
$$

Let us introduce a functional $L(\gamma (\cdot))$ on $F$ by the formula
$$
L(\gamma (\cdot))=\exp\Bigl(\int_X \log l(\gamma (x))\,d\,m(x)\Bigr),
$$
i.e.,
\begin{equation}{}\label{11}
L(\gamma (\cdot))=\exp\Bigl(-\frac12\int_X
\log\Bigl(1+\frac14 |\gamma (x)|^2\Bigr)\,d\,m(x)\Bigr).
\end{equation}

Note that
$
L(\gamma (\cdot))=\Psi(\gamma (\cdot)),
$
where $\Psi(\gamma (\cdot))$ is the restriction of the spherical function
of the basic representation of the group
$G^X$ to the subgroup $Z^X$, see \eqref{1-12}.
This functional is positive definite and continuous
(see \cite{V-G-7}); hence, by the Minlos theorem on measures on the
space of distributions
\cite{V-G-7}, it is the Fourier transform of a finite
normalized measure
$\mu$ on $\Phi$, i.e.,
\begin{equation}{}\label{12}
L(\gamma (\cdot)) =\int_{\textstyle \Phi} e^{\,i< \xi , \gamma >}\,d\mu(\xi ).
\end{equation}

\begin{THM}{}\label{THM:11} The measure $\mu$ is concentrated on the set
$\Phi_0 \subset \Phi$ of distributions of the form
\begin{equation}{}\label{727}
\xi=\sum c^i\, \delta _{x_i},\quad
c^i=(c^i_1, \dots ,c^i_{n-1})\in\RR^{n-1},\quad \text{where}
\quad  \sum |c^i|< \infty .
\end{equation}
\end{THM}
\begin{proof}{}
The series in \eqref{727} converges if and only if
the series for each coordinate converges;
hence it suffice to verify the 
condition for the one-dimensional
processes obtained by projecting to coordinates.
The characteristic
functions of these processes are the restrictions of the function
$l(\cdot)$ (given by \eqref{181-1}) to the one-dimensional subspaces; but 
all of them determine the classical gamma process, which satisfies the
convergence condition
(see \cite{V-8}; for a general convergence condition for one-dimensional 
L\'evy processes, see \cite{S-14}).
\end{proof}

Let us introduce the space 
$\Cal F \supset F $ of real bounded Borel vector functions on $X$.
Since $\mu$-almost every distribution
$\xi\in\Phi$ is of the form \eqref{727},
each function $\gamma \in\Cal F$ corresponds to a measurable linear functional
$<\xi, \gamma >$ on $\Phi$ defined $\mu$-almost everywhere on $\Phi$.

Note that on  $\Cal F$ and $\Phi$ there are two natural operations :

1) the multiplication by a bounded Borel  $\RR^*$-valued function
$\epsilon(x)$ on
$X$:
$$
\gamma \to \epsilon \gamma ,\quad  \xi\to \epsilon \xi;
$$

2) the action of the group $O(n{-}1)^X$ of Borel functions on $X$ with
values in the compact group $O(n{-}1)$:
$$
\gamma \to \gamma u,\quad \xi\to\xi u, \quad u\in O(n{-}1)^X.
$$

According to this definition,
\begin{align*}{}
&< \epsilon \xi, \gamma >\,=\,<  \xi, \epsilon\gamma >
&&\text{and}
&&< \xi u, \gamma >\,=\,<  \xi, \gamma u>
&&\text{for any}
&&\gamma \in\Cal F, \quad \xi\in\Phi.
\end{align*}

The following proposition follows from the definition of the measure $\mu$.

\begin{PROP}{}\label{PROP:213-5} The measure $\mu$ is invariant under the
transformations $\xi \to \xi u$, $u\in O(n{-}1)^X$.
\end{PROP}

\begin{REM*}{} 
In fact, we construct a vector gamma process such that the measure
in the space of trajectories of this process enjoys the additional property
as compared with the symmetrized gamma process: it is invariant
under pointwise orthogonal transformations.
\end{REM*}

\subsection{Projections of the measure $\mu$ 
to finite-dimensional quotient spaces}

Let us consider all finite partitions $X=\bigcup_{i=1}^l X_i$ of the space
$X$ such that $m(X_i)<n-1$ for all elements of the partition.
Let us associate with each partition
$\alpha:  X=\bigcup_{i=1}^l X_i$, where
$m(X_i)< n-1$, the subspace
$\Cal F_{\alpha }\cong (\RR^{n-1})^l$ of vector functions
$\gamma (x)\in\Cal F$ that are constant on elements of 
$\alpha$, and the dual space
$\Phi_{\alpha }\cong (\RR^{n-1})^l$, realized as a quotient space of 
$\Phi$.

Denote by $\mu_{\alpha }$  the projection of the measure $\mu$ 
to the quotient space
$\Phi_{\alpha }$.

\begin{PROP}{}\label{00-PROP:12}
The measure $\mu_{\alpha }$ on $\Phi_{\alpha }\cong (\RR^{n-1})^l$
has the form
\begin{equation}{}\label{17-9}
d\mu_{\alpha }(\xi_1, \dots ,\xi_l)=
\prod_{k=1}^l\Bigl(\frac{2}{\Gamma(\lambda_k /2)}|\xi_k|^{(\lambda_k -n+1)/2}\,
K_{\frac{n-1- \lambda_k}{2}}(2|\xi_k|)\,d\xi_k\Bigr),
\end{equation}
where $\xi_i\in\RR^{n-1}$, $d\xi_i$ is the Lebesgue measure on $\RR^{n-1}$,
$\lambda _k=m(X_k)$, and $K_{\nu}$
is the modified Bessel function of the third kind defined above.
\end{PROP}

\begin{proof}{}  For every $\gamma (x)\in\Cal F_{\alpha }$ we have
$$
\exp\Bigl(-\frac12\int_X \log\Bigl(1+\frac{|\gamma(x)|^2}{4}\Bigr)\,dm(x)\Bigl)
= \prod_{i=1}^l \Bigl(1+\frac{|\gamma^i|^2 }{4}\Bigr)^{- \lambda _i/2},
$$
where $\gamma^i= \gamma (x)|_{X_i}$ and $\lambda _i=m(X_i)$.

Therefore, in view of \eqref{12}, for every
$\gamma(x) \in\Cal F_{\alpha }$ we have
$$
\int_{\Phi_{\alpha }} e^{\,i<\xi, \gamma >} d\mu_{\alpha }(\xi)=
\int_{\Phi} e^{\,i<\xi, \gamma >} d\mu(\xi)=
\prod_{i=1}^l \Bigl(1+\frac{|\gamma^i|^2 }{4}\Bigr)^{- \lambda _i/2}.
$$
Hence $d\mu_{\alpha }(\xi^1, \dots ,\xi^l)=
\prod_{i=1}^l \psi(\xi^i)\,d\xi^i$, where $\psi(\xi^i)$ is the Fourier
transform of the function
$(1+\frac{|\gamma^i|^2 }{4})^{- \lambda _i/2}$ on $\RR^{n-1}$. Now
 \eqref{17-9} follows from \eqref{17-3}.
\end{proof}

\subsection{The function $V_{\rho}(x)$ and the measure 
$\nu$ on $\Phi$}

For $\rho>0$, 
introduce a function $V_{\rho }(x)$
on the half-line $0\le x< \infty $ by the formula
\begin{equation}{}\label{17-7}
V_{\rho}(x) = \Bigl(
        \frac{2}{\Gamma (\rho)}x^{\rho}\,K_{\rho}(2x)
        \Bigr)^{-1}
= \bigl(
        \Gamma (1-\rho)\,x^{\rho}\,[I_{-\rho}(2x)-I_{\rho}(2x)]
        \bigr)^{-1}.
\end{equation}
In particular, $V_{1/2}(x)=e^{2x}$.

\begin{THM}{}\label{THM:717-7} The Fourier transform of the function
$V_{\rho}^{-1}(|\xi|)$ on $\RR^{n-1}$ equals
\begin{equation}{}\label{727-7}
\int_{\RR^{n-1}} V_{\rho}^{-1}(|\xi|)\,e^{\,i\,<\xi,x>}\,d\xi=
c_n\,        \frac{\Gamma (\rho)}{\Gamma (\frac{n-1}{2}+\rho)}\,
\Bigl(1+\frac{|x|^2}{4}\Bigr)^{-\frac{n-1}{2}-\rho}.
\end{equation}
In particular,
$$
\int_{\RR^{n-1}} V_{\frac{n-1}{2}}^{-1}(|\xi|)\,e^{\,i\,<\xi,x>}\,d\xi=
c_n\,\frac{\Gamma (\frac{n-1}{2})}{\Gamma (n-1)}\,
\Bigl(1+\frac{|x|^2}{4}\Bigr)^{-n+1}.
$$
\end{THM}

Indeed, according to \eqref{717-3},
\begin{equation}{}\label{717-7}
V_{\rho}^{-1}(|\xi|) = c_n^{-1}\,
        \frac{\Gamma (\frac{n-1}{2}+\rho)}{\Gamma (\rho)}\,
\int_{\RR^{n-1}} \Bigl(1+\frac{|x|^2}{4}\Bigr)^{-\frac{n-1}{2}-\rho}\,e^{\,i\,<\xi,x>}\,dx.
\end{equation}
Applying the inverse Fourier transform yields \eqref{727-7}.

\begin{COR*}{} The measure
$V_{\frac{n-1}{2}}^{-1}(|\xi|)\,d\xi$ is an infinitely
divisible measure on
$\RR^{n-1}$.
\end{COR*}

(Since its Fourier transform has a L\'evy--Khintchin representation; see 
\eqref{345-1}.)

\begin{PROP}{}\label{PROP:456}
The function $V_{\rho }(x)$ is continuous and strictry positive on the half-line
$0\le x< \infty$, satisfies 
$V_{\rho} (0)=1$ for every $\rho>0$, and 
has the following asymptotic estimates as $x\to0$:
\begin{equation}{}\label{17-8}
V_{\rho}(x)\sim \begin{cases}{}
1+x^{\,2\rho}\,\Gamma (1-\rho)/\Gamma (1+\rho)&\text{for $\rho<1$},\\
1-2x^{\,2}\,\log(x)&\text{for $\rho=1$},\\
1+2x^{\,2}/(\rho-1)&\text{for $\rho>1$}.
\end{cases}
\end{equation}
\end{PROP}

\begin{proof}{} Let us prove \eqref{17-8}. If $\rho\notin\ZZ$, then
we use the estimate for the functions
$I_{-\rho}$ and $I_{\rho}$ that follows from their power series 
representation:
$$
I_{-\rho}(2x)\sim\frac{x^{-\rho}}{\Gamma (1-\rho)}
+\frac{x^{2-\rho}}{\Gamma (2-\rho)}\,x^2,\quad
I_{\rho}(2x)\sim\frac{x^{\rho}}{\Gamma (1+\rho)}.
$$
Hence we have
$$
V^{-1}_{\rho}(x)\sim 1+\frac{\Gamma (1-\rho)}{\Gamma (2-\rho)}\,x^2-
       \frac{ \Gamma (1-\rho) }{ \Gamma (1+\rho) }\,x^{\,2\rho}.
$$
This estimate implies \eqref{17-8} for $\rho\notin\ZZ$.

If $\rho\in\ZZ$, then we use the series representation of 
$K_n(2x)$, see
\cite[\S7.2.5, formula (37)]{B-12}. This representation
implies the estimate
$$
2\,K_n(2x)\sim (n-1)!x^{-n} +(n-2)!\,x^{2-n}+
\frac{2(-1)^{n+1}}{n!}x^n\log x.
$$
Therefore
$$
V^{-1}_n(x)\sim 1+\frac{x^2}{n-1}+\frac{2(-1)^{n+1}}{n!(n-1)!}x^{2n}\log x,
$$
whence $V^{-1}_1\sim 1+2 x^2\log x$ and
$V^{-1}_n\sim 1+\frac{x^2}{n-1}$ for $n>1$.
\end{proof}

\begin{COR*}{} For every $\rho\ge 1/2$, the infinite product
$ \prod_{i=1}^ \infty V_{\rho}(x_i)$ converges provided that the series 
$ \sum_{i=1}^ \infty x_i$ converges.
\end{COR*}

\begin{DEF*}{} Consider the infinite $\sigma $-finite measure $\nu$ on
the space of vector distributions 
$\Phi$ whose density
$v=\frac{d\nu}{d\mu}$ with respect to the measure $\mu$ is defined on
the support
$\Phi_0 \subset \Phi$ of $\mu$ by the following formula:
\begin{equation}{}\label{29}
v\bigl(\sum_{i=1}^\infty c^i \delta _{x_i}\bigr)=2^{-m(X)}\prod_{i=1}^\infty
V_{\frac{n-1}{2}}(|c^i|), \quad c^i\in\RR^{n-1},
\end{equation}
where $V_{\rho}(x)$ is given by \eqref{17-7}.
In particular, for $n=2$
$$
v\bigl(\sum c^i \delta _{x_i}\bigr)=\exp(2\sum |c^i|),
\quad c^i\in\RR.
$$
\end{DEF*}

Since $\sum_{i=1}^\infty |c^i| < \infty $ on the support of $\mu$, 
it follows from Proposition~\ref{PROP:456} that the infinite product
$\prod_{i=1}^\infty V_{\frac{n-1}{2}}(|c^i|)$ converges.

By definition, the measure $\nu$ 
is absolutely continuous with respect to
$\mu$, and its density is positive $\mu$-almost everywhere.

\subsection{Approximative construction of the measure $\nu$}
Let $\alpha:  X=\bigcup_{i=1}^l X_i $ be an arbitrary finite partition of
the space $X$, $\Phi_{\alpha }$ be the quotient space of
$\Phi$ associated with $\alpha$, and
$\mu_{\alpha }$ be the projection of the measure $\mu$ to $\Phi_{\alpha }$.
Let us introduce a new measure $\nu_{\alpha }$  on
$\Phi_{\alpha }$ with density
$$
\frac{d\nu_{\alpha }(\xi^1, \dots ,\xi^l)}
{d\mu_{\alpha }(\xi^1, \dots ,\xi^l)}=2^{-m(X)}
\prod_{k=1}^l V_{(n-1- \lambda _k)/2}(|\xi^k|).
$$
(By the condition imposed on $\alpha$, we have
$n-1- \lambda _k>0$ for all $k$.)

It follows from the explicit formula for $d\mu_{\alpha }(\xi^1, \dots ,\xi^l)$ 
that this measure is of the form
\begin{equation}{}\label{17-91}
d\nu_{\alpha }(\xi^1, \dots ,\xi^l)=
\prod_{i=1}^l \frac{2^{\,-\lambda_i}\, \Gamma ((n-1-\lambda_i)/2) }
{\Gamma (\lambda_i /2)}\,|\xi^i|^{\lambda_i -n+1}\,d\xi^i,
\end{equation}
where $\xi^i\in\RR^{n-1}$ and
$\lambda _i=m(X_i)$.

Let us write $\beta \ge \alpha$ if $\beta $ is a refinement of a partition
$\alpha$.  For $\beta\ge \alpha $, there is a natural embedding
$\Cal F_{\alpha }\to\Cal F_{\beta }$ and a natural epimorphism
$\Phi_{\beta }\to\Phi_{\alpha }$. Obviously, the measure
$\mu_{\beta }$ on
$\Phi_{\beta }$ and the measure
$\mu_{\alpha }$ on $\Phi_{\alpha }$ 
are coherent with respect to this epimorphism.

\begin{PROP}{}\label{PROP:29} The measures $\nu_{\alpha }$ 
on the quotient spaces
$\Phi_{\alpha }$ are coherent, i.e., the epimorphism
$\Phi _{\beta }\to \Phi _{\alpha }$, $\beta \ge \alpha $, sends 
$\nu _{\beta }$ to $\nu _{\alpha }$.
\end{PROP}
\begin{proof}{} It suffices to prove that
$$
\int_{\Phi_{\alpha }} e^{\,i<\xi, \gamma >}\,d\nu_{\alpha }(\xi)=
\int_{\Phi_{\beta }} e^{\,i<\xi, \gamma >}\,d\nu_{\beta }(\xi) \quad
\text{for every} \quad \gamma \in\Cal F_{\alpha }.
$$
Let $\alpha : X=\bigcup_{i=1}^lX_{i}$ and
$\beta :  X=\bigcup_{i,j}X_{ij}$, where 
$\bigcup_{j} X_{ij}=X_i,$ $i=1, \dots ,l$. By \eqref{17-4}, it follows
from the expression \eqref{17-91} for $d \nu _{\alpha }$ and $d \mu _{\alpha }$
that
\begin{eqnarray*}{}
\int_{\Phi_{\alpha }} e^{\,i<\xi, \gamma >}\,d\nu_{\alpha }(\xi)&=&
\prod_{i=1}^l |\gamma^i|^{- \lambda _i}, \quad \text{where}\quad
\gamma^i= \gamma |_{X_i}, \quad \lambda _i=m(X_i);
\\
\int_{\Phi_{\beta }} e^{\,i<\xi, \gamma >}\,d\nu_{\beta }(\xi)&=&
\prod_{i,j} |\gamma ^{\,ij}|^{- \lambda _{ij}}, \quad \text{where}\quad
\gamma ^{\,ij}= \gamma |_{X_{ij}}, \quad \lambda _{ij}=m(X_{ij}).
\end{eqnarray*}
Since $\gamma ^{\,ij}= \gamma^i$ for all $i$ and $j$ and
$\sum_j \lambda_{ij}= \lambda _i,$ $i=1, \dots ,l$,
the right-hand sides of these equations coincide.
\end{proof}

\begin{DEF*}{} Let us define a $\sigma $-finite measure $\TIL\nu$ on $\Phi$
as the weak limit of the coherent family of measures
$\nu_{\alpha }$.
\end{DEF*}

\begin{THM}{}\label{THM:14}  The measure $\TIL\nu$ coincides with the measure
$\nu$, i.e., on the support of $\mu$,
\begin{equation}{}\label{299}
\frac{d\TIL\nu}{d\mu}\bigl(\sum c^i \delta _{x_i}\bigr)=2^{-m(X)}
\prod V_{\frac{n-1}{2}}(|c^i|),
\quad c^i\in\RR^{n-1},
\end{equation}
where
$V_{\rho}(x)$ is given by \eqref{17-7}.
\end{THM}

\begin{proof}{}
It suffices to prove \eqref{299} only for finite sums
$\xi=\sum_{i=1}^k c^i\, \delta _{x_i}$.

Given such a sum, consider partitions
$\alpha :  X=\bigcup_{i=1}^l X_i$ such that each element of $\alpha$
contains at most one point $x_i$. For definiteness, let
$x_i\in X_i$, $i=1, \dots ,k$. Then, since $V_{\rho}(0)=1$, we have
$$
\frac{d\nu_{\alpha }}{d\mu_{\alpha }}\Bigl(
\sum_{i=1}^k c^i \delta _{x_i}\Bigr)=2^{-m(X)}\,
\prod_{i=1}^k V_{\frac{n-1-\lambda_i}{2}}(|c^i|).
$$
Taking the inductive limit with respect to $\alpha$, we obtain
$$
\frac{d\TIL\nu}{d\mu}\Bigl(
\sum_{i=1}^k c^i \delta _{x_i}\Bigr)=2^{-m(X)}\,
\prod_{i=1}^k V_{\frac{n-1}{2}}(|c^i|).
$$
\end{proof}

\subsection{The Fourier transform of the measure $\nu$}
If $\gamma \in \Cal F_{\alpha }$, where $\alpha $ 
is an arbitrary finite partition of the space $X$, then by
\eqref{17-4} we have
$$
\int_{\Phi} e^{\,i\,<\xi,\gamma>}\,d\nu(\xi)=
\int_{\Phi_{\alpha }} e^{\,i\,<\xi,\gamma>}\,d\nu_{\alpha }(\xi)=
\prod_{i=1}^l |\gamma^i|^{- \lambda _i},
$$
where $\gamma^i= \gamma |_{X_i}$ and $\lambda _i=m(X_i)$.
The right-hand side of this equation can be represented in the form
$$
\prod_{i=1}^l |\gamma^i|^{- \lambda _i}=
\exp\Bigl(-\int_X\log |\gamma(x)|\,dm(x)\Bigr).
$$
Thus
\begin{equation}{}\label{30-1}
\int_{\Phi} e^{\,i\,<\xi,\gamma>}\,d\nu(\xi)=
\exp\Bigl(-\int_X\log |\gamma(x)|\,dm(x)\Bigr).
\end{equation}
Equation \eqref{30-1} determines the Fourier transform of the measure
$\nu$ whenever the integral in the right-hand side is finite.

\begin{REM*}{} One can take this equation as the definition of the measure
$\nu$.
\end{REM*}

\subsection{The invariance properties of the measure $\nu$}

\begin{THM}{}\label{THM:44}
The measure $\nu$ is invariant under the action of the group
$O(n{-}1)^X$ and projective invariant under the multiplication by
bounded Borel functions
$\epsilon (x)\in(\RR^*)^X$  such that
the integral $\int_X \log |\epsilon (x)|\,d m(x)$ converges. Namely,
\begin{eqnarray}{}\label{18-33}
d\nu(\xi\,u)&=&d\nu(\xi) \quad
\text{for every} \quad u(x)\in O(n{-}1)^X;
\\
\label{18-3}
d\nu( \epsilon\,\xi )&=&e^{\, \int_X \log |\epsilon(x)|\,dm(x)}\,d\nu(\xi).
\end{eqnarray}
In particular, $\nu$ is invariant under the subgroup of 
multiplications by functions
$\epsilon (x)$ satysfying 
$\int_X \log |\epsilon(x)|\,dm(x)=0$.
\end{THM}

\begin{proof}{} It suffices to establish this property for the projections
$d\nu_{\alpha }(\xi)=d\nu_{\alpha }(\xi_1, \dots ,\xi_l)$ of $\nu$
to the quotient spaces $\Phi_{\alpha }$ of the space $\Phi$. Let
$u(x)\in O(n{-}1)^X$ and $\epsilon (x)\in(\RR^*)^X$ 
be constant on the elements of a partition
$\alpha$. By definition,
\begin{gather*}{}
d\nu_{\alpha }(\xi u)
=d\nu_{\alpha }(\xi_1u_1, \dots ,\xi_lu_l),
\\
d\nu_{\alpha }(\epsilon  \xi)
=d\nu_{\alpha }(\epsilon_1 \xi_1, \dots ,\epsilon_l \xi_l),
\end{gather*}
where $u_i=u|_{X_i}$, $\epsilon _i= \epsilon |_{X_i}$, $X_i$ are elements
of $\alpha$.

It follows immediately from the explicit expression \eqref{17-91} for
$d\nu_{\alpha }$ that
$d\nu_{\alpha }(\xi u)=d\nu_{\alpha }(\xi)$ and
$d\nu_{\alpha }(\epsilon  \xi)=\prod_{i=1}^l |\xi_i|^{\lambda _i}\,
d\nu_{\alpha }(\epsilon  \xi)$, where $\lambda _i=m(X_i)$. It suffices
to observe that
$
\displaystyle
\prod_{i=1}^l |\xi_i|^{\lambda _i}=
\int_X \log |\epsilon(x)|\,dm(x)
\,d\nu(\xi).
$
\end{proof}

\section{Construction of the basic representation of the current group
$G^X$}

\subsection{The basic representation of the 
block triangular group $B^X$} 

First let us describe the representation of the subgroup
$B^X=Z^X\leftthreetimes D^X$ of block triangular matrices.
We will write elements of this subgroup as triples
$$
g(x)=(\epsilon ,u, \gamma ),\quad \epsilon \in(\RR^*)^X ,
\quad u(x)\in O(n{-}1)^X,\quad \gamma \in (\RR^{n-1})^X.
$$
In this notation, the product of group elements takes the form
$$
(\epsilon_1 ,u_1, \gamma_1)\,(\epsilon_2 ,u_2, \gamma_2)=
(\epsilon_1 \epsilon _2,\, u_1u_2,\, \gamma_1+ \epsilon _1 \gamma _2u_1^{-1}).
$$

The representation of the group $B^X$ is realized in the Hilbert space
$L^2(\nu)$  of all functions on $\Phi$ square-integrable with respect to
the measure $\nu$ introduced above.

Let us associate with elements of the subgroup
$B^X$ the following operators $U_g$ in the space of functions
$f(\xi)$ on $\Phi$:
\begin{equation}{}\label{17-13}
U_{\epsilon ,u,\gamma} f(\xi)=e^{\,1/2\int_X\log|\epsilon (x)|\,dm(x)+
i\,<\xi,\gamma>}\,f(\epsilon\,\xi u).
\end{equation}
In particular, elements $z=z(\gamma )\in Z^X$ give rise to the
operators
$$
U_z f(\xi)=e^{\,i\,<\xi,\gamma>}\,f(\xi ),
$$
where $\gamma\in(\RR^{n-1})^X $ is the parameter of
the block matrix $z$, 
and elements $d=\diag(\epsilon ^{-1}, u, \epsilon )\in D^X$, to 
the operators
$$
U_d f(\xi)=e^{\,1/2\int_X\log|\epsilon (x)|\,dm(x)}
\,f(\epsilon\,\xi u).
$$

Note that the integral $\int_X \log |\epsilon (x)|\,dm(x)$
converges, because the elements
$\epsilon (x)$ and $\epsilon^{-1} (x)$
of the matrix $d$ are bounded functions on $X$.

\begin{THM}{}\label{THM:2}
The operators $U_g$, $g \in B^X$, given by \eqref{17-13}
are unitary with respect to the norm in
$L^2(\nu)$ and form an irreducible representation of the group
$B^X$ in $L^2(\nu )$.
\end{THM}

Indeed, the group property of the operators
$U_g$ follows immediately from their definition, and the unitarity
follows from the invariance properties of the measure
$\nu$ established in Theorem~\ref{THM:44}. The representation
$U_g$ is irreducible, because the action of the subgroup
$Z^X$ in $L^2(\nu )$ is ergodic and the algebra of multiplicators is maximal.

\subsection{Extension of the representation $U_g$ to
the whole group $G^X$}

In order to extend the representation
$U_g$ from the subgroup $B^X$ to the whole group $G^X$,
we use the following obvious lemma.

\begin{LEM*}{}The group $O(n,1)^X$ is algebraically generated by the subgroup
$B^X$ and the unique element
$g(x)\equiv s$.
\end{LEM*}

For example, every element $g$ of the block upper triangular subgroup
can be represented in the form
$g=sg^*s$, where $g^*\in B^X$.

By this lemma, every operator $U_g$, $g\in G^X$, can be represented as
the product of an operator from the subgroup
$B^X$ described above and the operator
$U_s$; hence, in order to define the representation of the whole
group $G^X$, it suffices to describe only the operator $U_s$.

We call $U_s$ the involution operator in the space
$L^2(\nu)$ and denote it by $I$.

\begin{PROP}{}\label{PROP:98}
The operator $I=U_s$ and
the operators $U_g$, $g\in B^X$, satisfy
the following relations:
\begin{eqnarray}
\label{363}
I\,U_d&=&U_{s\,d\,s}\,I \quad \text{for every} \quad d\in D^X,
\\
\label{364}
U_{z(\gamma )}\,I&=&U_{d(\gamma )}\,I\,U_{z(-\gamma )}\,I\,U_{z(j\gamma )}
\quad \text{for every} \quad z(\gamma )\in Z^X,
\end{eqnarray}
where $j\gamma =-\frac{2 \gamma }{|\gamma |^2}$ (involution) and
$d(\gamma )\in D^X$ is given by 
\begin{equation}{}\label{391}
d(\gamma )=\diag\Bigl(-\frac{2}{|\gamma |^2 },\,
u_{\,\gamma },\,-\frac{|\gamma |^2 }{2}\Bigr),
\qquad
u_{\gamma }=e-\frac{2 \gamma ^* \gamma }{|\gamma |^2 }.
\end{equation}
\end{PROP}

Indeed,  \eqref{363} is obvious and \eqref{364} follows from
the corresponding relation in the group $O(n,1)^X$:
$$
z(\gamma )\,s=d(\gamma )\,s\,z(- \gamma )\,s\,z(j \gamma ).
$$

Note that for every $x\in X$ the matrix $u_{\gamma }(x)$
determines the reflection in
$\RR^{n-1}$ with respect to the hyperplane orthogonal to the vector
$\gamma (x)$.

Relations \eqref{363} and \eqref{364} uniquely determine the operator
$I$; however, they do not give an explicit expression for this operator.

An explicit descriprion for $I$
can be obtained from the description of the operator
$T_s$ in the commutative model of complementary series
representations of the group
$O(n,1)$. Namely, let us associate with each partition 
$\alpha : X=\bigcup_{i=1}^l X_i$ the Hilbert space
$L^2(\nu_{\alpha }) \subset L^2(\nu)$
of functions on $\Phi_{\alpha }$ square-integrable with respect to the measure
$\nu_{a}$. The spaces
$L^2(\nu_{\alpha })$ are invariant under $I$, 
and in order to describe the operator
$I$ on the whole space $L^2(\nu)$, it suffices to describe its action
on each of these subspaces. By construction, each space
$L^2(\nu_{\alpha })$ is isomorphic to the tensor product of finitely many spaces
on which the commutative model of a complementary series
representation of the group $O(n,1)$ acts.

\begin{PROP}{}\label{PROP:578} On each subspace $L^2(\nu_{\alpha })$, 
the involution operator
$I$ coincides with the operator $T_s$ of the representation of the group
$O(n,1)$ in this subspace.
\end{PROP}

Starting from formulas \eqref{73-1} and \eqref{63-1} for the operator $T_s$ in
the commutative model of complementary series representations of the group
$O(n,1)$, we obtain the following theorem.

\begin{THM}{}\label{THM:78} For every partition
$\alpha : X=\bigcup_{i=1}^l X_i$, the action of the operator $I$ on the subspace
$L^2(\nu_{\alpha })$ is given by the following formula:
$$
If(\xi_1, \dots ,\xi_l)
=\!\!\!\!\!\!\int\limits_{(\RR^{n-1})^l}\!\!\!\!\!\!
\Bigl(f(\xi'_1, \dots ,\xi'_l)\prod_{i=1}^l
A^{\lambda _i} (\xi_i,\,\xi'_i)\,d\xi_i\Bigr),\quad \xi_i\in\RR^{n-1}, \quad
\lambda _i=m(X_i),
$$
where
$$
A^{\lambda }(\xi,\xi')=c_n\,2^{- \lambda/2 }\,
\int_0^\infty r^{\,\lambda -n}\,|r\xi+2\,r^{-1}\xi'|^{-\frac{n-3}{2}}\,
J_{\frac{n-3}{2}}(|r\xi+2\,r^{-1}\xi'|)\,dr.
$$
\end{THM}

We see that the operator $I$ is well defined on the whole space
$L^2(\nu)$, because on all subspaces
$L^2(\nu_{\alpha})$ it is unitary
and satisfies the required relations.
Thus we have constructed
an irreducible unitary representation of the group
$G^X$ in the space
$L^2(\nu)$.

\begin{THM}{}\label{THM:22} The constructed representation 
$ U_g$ of the group $G^X$ is equivalent to the basic representation
of this group introduced in
\cite{V-G-6}.
\end{THM}

\begin{proof}{} The representations are equivalent, because 
their spherical functions coincide on the subgroup
$Z^X$. Namely, set
$\phi(\xi)=v^{-1/2}(\xi)$, where $v(\xi)$ is the density of the measure $\nu$
with respect to $\mu$. The function $\phi$ belongs to the space $L^2(\nu)$,
is of norm $1$,
and satisfies the equation
$$
< U_{z(\gamma )\,} v^{-1/2},\, v^{-1/2}>\,=\Psi(\gamma (\cdot)),
\quad z(\gamma )\in Z^X,
$$
where $\Psi(\gamma (\cdot))$ is the restriction of the spherical function
of the basic representation of
$G^X$ to $Z^X$. Indeed,
$$
< U_{z(\gamma )} v^{-1/2},\, v^{-1/2}>\,=\int_{\Phi} e^{\,i<\xi, \gamma >}\,
v^{-1}(\xi)\,d\nu(\xi)=\int_{\Phi} e^{\,i<\xi, \gamma >}\,d\mu(\xi)=
\Psi(\gamma (\cdot)).
$$
\end{proof}

\begin{REM*}{} One can also construct a commutative model of the basic 
representation of the group $G^X$ in the 
$L^2$ space over the probability measure $ \mu$. However, 
$\mu$ is only quasi-invariant under the transformations 
$\xi\to \epsilon \xi$, $\epsilon \in (\RR^*)^X$.
Thus, in order to obtain a unitary representation, one must
introduce an additional factor. In the case of the $\sigma$-finite 
measure $\nu$, there is no need to to this.
\end{REM*}

\section{Appendices}

\subsection{Approximative construction of the commutative model}

Let us give another, independent  construction of the commutative model of the 
basic representation of the group 
$G^X$. With each finite 
partition $\alpha:  X=\bigcup_{i=1}^l X_i $ of the space $X$
associate the Hilbert space $L^2(\nu_{\alpha })$ of functions on 
$\Phi_{\alpha }$ introduced above and the subgroup
$G^X_{\alpha }\subset G^X$ of functions constant on the elements
of $\alpha$.

For $\beta \ge \alpha $, there is a natural embedding of groups
$G^X_{\alpha } \to  G^X_{\beta }$ and a natural isometric embedding of spaces
$L^2(\nu_{\alpha }) \to L^2(\nu_{\beta })$.
Denote by $L_0$ the inductive limit of the subspaces
$L^2(\nu_{\alpha })$, and by $G_0^X$ the inductive limit of the subgroups
$G^X_{\alpha }$. Note that the group $G_0^X$ is everywhere
dense in $G^X$.

We will define a unitary representation $U_g$ 
of the group $G^X$ on the completion
$L$ of the space
$L_0$ with respect to the norm of $L_0$.  In order to define it, it suffices
to describe the action of the operators
$U_g$, $g\in G^X_0$, on the subspace $L_0$.

\begin{PROP}{}\label{PROP:28} The space $L^2(\nu_{\alpha})$
is the tensor product
\begin{equation}{}\label{4}
L^2(\nu_{\alpha })= \bigotimes_{i=1}^l L_{\lambda _i},
\quad \lambda _i=m(X_i),
\end{equation}
where $L_{\lambda _i}$ are the Hilbert spaces introduced in
Theorem~{\rm \ref{THM:1}}.
\end{PROP}

Starting from the decomposition \eqref{4}, define a unitary
representation $U_g$ of the group $G^X_{\alpha }$ in $L^2(\nu_{a})$ 
by the formula
\begin{equation}{}\label{63-2}
U_{g(x)}=T^{\lambda _1}_{g_1}\otimes \ldots
\otimes T^{\lambda _l}_{g_l}, \quad g_i=g(x)|_{X_i},
\end{equation}
where $T^{\lambda _i}$ are operators of the complementary series representations
of the group $G=O(n,1)$ in the spaces $L_{\lambda _i}$ defined
in Theorem~\ref{THM:1}.

\begin{PROP}{}\label{PROP:71} For $\beta \ge \alpha $,
the embedding $L^2(\nu_{\alpha }) \to L^2(\nu_{\beta })$
commutes with the representations $\TIL U_g$ of the groups 
$G^X_{\alpha }$ and $G^X_{\beta }$
in these spaces.
\end{PROP}

\begin{COR*}{} The unitary representations $U_g$ of the groups $G^X_{\alpha }$ in
the spaces $L^2(\nu_{\alpha })$ are coherent and hence generate a unitary
representation of the group
$G_0^X$ in $L_0$.
\end{COR*}

\begin{THM}{}\label{THM:4} The contructed representation $U_g$ of the group $G^X$
is irreducible and equivalent to the basic representation of this group
defined in \cite{V-G-2,V-G-6}.
\end{THM}

\subsection{The dual description of the representation $U_g$} 
Let us give a description of the representation
$U_g$ in terms of the Fourier transform sending functions
$f(\xi)$
on $\Phi$ to functions $\phi(\gamma )$ on $\Cal F$:
$$
\Cal R f(\gamma )=
\int_{\textstyle \Phi} f(\xi)\,e^{\,i<\xi, \gamma >}\,d\nu(\xi),\quad
\gamma \in\Cal F.
$$
The operator $\Cal R$ is defined on an everywhere dense subset of functions
$f\in L^2(\nu)$; in particular, if the integral
$\int_{\Phi} |f(\xi)|\,d\nu(\xi)$ converges, then the function
$\Cal R f(\gamma )$ is defined on the whole space
$\Cal F$.

\begin{THM}{}\label{THM:25}
Let $f_1\in L^2(\nu)$, $g \in G^X$, and $f_2=U_g f_1$.
Then on the subset of $\gamma \in\Cal F$ for which the functions $\Cal R f_1$ and
$\Cal R f_2$ are defined, they satisfy the relation
\begin{equation}{}\label{257}
\Cal R f_2(\gamma )=\Cal R f_1(\gamma\OVER g )\,e^{-1/2\,\int_X
\log \beta (\gamma (x),g(x))\,d m(x)},
\end{equation}
where, according to formulas \eqref{1-2} and \eqref{1-3}
for the group $G=O(n,1)$,
\begin{multline*}{}
\gamma(x)\,\OVER{g(x)}
\\
=\Bigl(-\frac{|\gamma(x)|^2}{2}\,g_{13}(x)+
\gamma(x) g_{23}(x)+g_{33}(x)\Bigr)^{-1}
\Bigl(-\frac{|\gamma(x)|^2}{2}\,g_{12}(x)+ \gamma(x) g_{22}(x)+g_{32}(x)\Bigr),
\end{multline*}
$$
\beta(\gamma (x),g(x))=\Bigl|-\frac{|\gamma(x)|^2}{2}g_{12}(x)+
\gamma(x)\, g_{22}(x)+g_{32}(x)\Bigr|.
$$
In particular,
\begin{align}{}
\label{444}
\Cal R f_2(\gamma )&=
\Cal R f_1(\gamma+ \gamma _0) \quad \text{for}\quad g=z(\gamma _0)\in Z^X;
\\
\notag
\Cal R f_2(\gamma )&=\Cal R f_1(\epsilon ^{-1} \gamma u)
\,e^{-1/2\,\int_X \log |\epsilon (x)|\,d m(x)}  \quad \text{for}\quad
g=\diag(\epsilon^{-1} ,u, \epsilon  )\in D^X;
\\
\notag
\Cal R f_2(\gamma )&=2^{\,1/2m(X)}\,\Cal R f_1\bigl(-\frac{2 \gamma }{|\gamma |^2}\bigr)
\,e^{-\int_X \log |\gamma (x)|\,d m(x)}  \quad \text{for}\quad g(x)\equiv s.
\end{align}
\end{THM}

\begin{proof}{} It suffices to prove \eqref{257} for functions
$f_1(\xi)=f_1(\xi^1, \dots ,\xi^l)$ from $L^2(\nu_{\alpha })$,
where $\alpha : \bigcup_{i=1}^l X_i$ is an arbitrary finite partition 
of the space $X$, and matrices
$g(x)$ that are constant on the elements of $\alpha$.

We have
$$
\Cal R f_i(\gamma )=\Cal R f_i(\gamma ^1, \dots , \gamma ^l)=
\int_{\Phi_{\alpha }}f_i(\xi^1, \dots ,\xi^l)\,e^{\,i<\xi, \gamma >}\,
d\nu_{\alpha }(\xi^1, \dots ,\xi^l),\quad i=1,2,
$$
where the measure $d\nu_{\alpha } (\xi^1, \dots ,\xi^l)$ is given by
\eqref{17-91}.  Since
$$
\int_{\Phi_{\alpha }} e^{\,i<\xi, \gamma >}\,
d\nu(\xi^1, \dots ,\xi^l)=\prod_{i=1}^l | \gamma^i|^{- \lambda _i},\quad
\lambda _i=m(X_i),
$$
it follows that
\begin{equation}\label{397}
\Cal R f_i(\gamma ^1, \dots , \gamma ^l)=
\int_{\Cal F_{\alpha }}\phi_i(\zeta^1, \dots , \zeta ^l)
\prod_{i=1}^l | \gamma^i- \zeta ^i|^{- \lambda _i}\,d\zeta^i,\quad i=1,2,
\end{equation}
where
$$
\phi_i(\zeta ^1, \dots , \zeta ^l)=\int_{\Cal F_{\alpha }}
\prod_{i=1}^l e^{\,i<\xi^i,\zeta^i>}\,f_i(\gamma ^1, \dots ,
\gamma ^l )\,d\zeta^i.
$$
The equation $f_2= U_g f_1$ and formula \eqref{1-8} for operators
of the complementary series representations of 
$O(n,1)$ imply that the functions
$\phi_i$ satisfy the relations
$$
\phi_2(\zeta ^1, \dots , \zeta ^l)=
\phi_1(\zeta^1\bar g_1, \dots , \zeta^l\bar g_l)
\,\prod_{i=1}^l \beta^{1-n+\frac{\lambda _i}{2}} (\gamma^i,g_i) ,
$$
where $g_i=g|_{X_i}$. Thus
$$
\Cal R f_2(\gamma ^1, \dots , \gamma ^l)=
\int_{\Cal F_{\alpha }}
\phi_1(\zeta^1\bar g_1, \dots , \zeta^l\bar g_l)
\,\prod_{i=1}^l \beta^{1-n+\frac{\lambda _i}{2}} (\gamma^i,g_i)\,
| \gamma^i- \zeta ^i|^{- \lambda _i}\,d\zeta^i.
$$
Applying the transformation $\zeta _i\bar g_i\to \zeta _i$, we obtain
$$
\Cal R f_2(\gamma ^1, \dots , \gamma ^l)=\prod_{i=1}^l
\beta ^{- \lambda _i/2}(\gamma^i,g_i)=
\int_{\Cal F_{\alpha }}\phi_1(\zeta^1, \dots , \zeta ^l)
\prod_{i=1}^l | \gamma^i\bar g_i- \zeta ^i|^{- \lambda _i}\,d\zeta^i,
$$
i.e.,
$$
\Cal R f_2(\gamma ^1, \dots , \gamma ^l)=
\Bigl(
\prod_{i=1}^l \beta ^{- \lambda _i/2}(\gamma^i,g_i)
\Bigr)
\Cal R f_1(\gamma ^1\bar g_1, \dots , \gamma ^l\bar g_l).
$$
It remains to observe that
$$
\prod_{i=1}^l \beta ^{- \lambda _i/2}(\gamma^i,g_i)=
e^{-1/2\,\int_X \log \beta (\gamma (x),g(x))\,d m(x)}.
$$
\end{proof}

\begin{COR*}{}
On the set of all $\gamma \in\Cal F$ for which the functions
$\Cal R f$, $f\in L^2(\nu)$, and $\Cal R I f$ are defined,
they satisfy the relation
\begin{equation}{}\label{2561}
\Cal R I f(\gamma )=2^{\,1/2m(X)}\,\Cal R f\Bigl(-\frac{2 \gamma }{|\gamma |^2}\Bigr)
\,e^{\,-\int_X \log |\gamma (x)|\,d m(x)}  \quad \text{for}\quad g(x)\equiv s.
\end{equation}
\end{COR*}

\subsection{On the properties of the measures $\mu $ and $\nu $}

\subsubsection{}
Consider the subgroup $O(n{-}1)^X$ of the current group $O(n,1)^X$. For  $n=2$,
this is the group of functions on $X$ taking values
$+1$ and $-1$.

The subgroup $O(n{-}1)^X$ acts pointwise in the space of vector 
distributions of dimension
$n{-}1$ on the manifold $X$, and the measures $\mu $ and $\nu $
are invariant under this action.

Note that this action is not free. Indeed, since almost every, with respect
to $\mu$ and $\nu$, realization is a linear combination of a countable
family of $\delta$-measures, it follows that for every distribution
$\xi $, only the compact quotient group of
$O(n{-}1)^X$ that consists of the restrictions of currents
to the countable support of $\xi$ acts freely on $\xi$.
Therefore the orbit is compact as the product of a countable family of
($n{-}1$)-dimensional spheres and hence has an invariant (product) measure.

Thus almost every ergodic component of the action of 
$O(n{-}1)^X$ consists of distributions with equal values of the norm
at all points of their common support, i.e., the function
$$
x \to \|c(x)\|,\quad c(x)=(c_1(x), \dots ,c_{n-1}(x)),
$$
is an invariant of the orbit.

At the same time, as noted above,
the action of the group of homotheties  
together with rotations, i.e., the group
$(\RR_+ \times O(n{-}1))^X$, on the space
$\Phi $ is already ergodic.
Recall that, as was proved in \cite{V-8} for the case $n=2$,
$\nu $ is the unique, up to normalization, measure that is
invariant and ergodic (Theorem~5).
Apparently, a similar theorem holds for an arbitrary $n$.

\subsubsection{}
The results of this paper can be translated to the infinite-dimensional
group $O(\infty ,1)^X$, because (as observed by
G.~Olshanski) all our products involve the dimension $n$ in a controllable way.
In other words, for different $n$ only
the natural dimension $n-1$ of the space of vector distribution changes,
but the form of the restriction of the spherical function to the subgroup
$Z^X$, and hence the characteristic functionals of the measures
$\mu $ and $\nu $ remain the same. Therefore our theory
can be considered for infinite $n$.

In this case, $\mu$-almost every vector distribution is still
a countable linear combination of $\delta$-measures on $X$,
but taking values in an infinite-dimensional Hilbert space
equipped with a mixture of Gaussian measures with some weight.

The role of the group $O(n{-}1)^X$ is played by the group $O(\infty)^X$
of orthogonal matrices of the form
$I+K$, where $K$ is a finite-dimensional operator. Remarkably, the
decomposition into ergodic components under the action of the group
$O(\infty)^X$ is the decomposition into the Gaussian measures with 
characteristic functionals of the form
$e^{-c\,\|x\|^2}$,
where $c$ is distributed on $(0,\infty )$ 
according to the measure with density $e^{-c}$.

\subsection{Remark on the group $U(n,1)^X$}
The above construction of an irreducible  unitary representation of the subgroup 
$B^X$ of $O(n,1)^X$ in the space
$L^2(\mu)$ can also be used for
the group
$U(n,1)^X$. However, in this case there appears a new phenomenon.

Like $O(n,1)$, we realize $U(n,1)$ as the group of linear transformations
in $\CC^{n+1}$ preserving the Hermitian form
$$
 x_1\OVER{x_{n+1}} + x_2\OVER{x_2} + \ldots + x_n\OVER{x_n}
$$
and represent its elements as block matrices. In this realization,
$B^X \subset U(n,1)^X$ is the semidirect product
$B^X=Z^X\leftthreetimes D^X$, where
$Z^X$ is the group of matrices of the form
$$
z=\begin{pmatrix}1&0&0\\
-\gamma^*&e&0\\it-\frac{|\gamma|^2}{2}&\gamma&1\end{pmatrix},
\quad
t\in\RR^X, \quad  \gamma \in (\CC^{n-1})^X
$$
(the Heisenberg group) and $D^X$
is the subgroup of block diagonal matrices 
$$
d=\diag(\bar\epsilon ^{-1}, u, \epsilon ), \quad \epsilon \in (\CC^*)^X,
\quad u\in U(n-1)^X.
$$
Accordingly, elements of the group
$B^X$ are 4-tuples
$b=(t, \gamma , \epsilon ,u)$.
If we now try to use a formular similar to 
\eqref{17-13},
$$
U_b \,f(\xi )
= e^{\,1/2\int_X\log|\epsilon (x)|\,dm(x) \,+\, i\,\opn{Re}<\xi,\gamma>}
\,f(\epsilon\,\xi u),
\qquad
b=(t, \gamma , \epsilon ,u),
$$
for representation operators in the space $L^2(\nu)$,
where the measure
$\nu$ is constructed as above, then we will find out that this
representation is not a {\it faithful representation}
of the group $B^X$,  and it is faithful only on the quotient
group with respect to the center of the Heisenberg subgroup,
since the operators of the representation do not involve the parameter
$t$. Thus the direct translation
of the construction does not use the simplectic 
structure on
$\CC^{n-1}\times \CC^{n-1}$, which is used in the definition of the group
$U(n,1)$, and cannot be extended to the whole group
$U(n,1)^X$. Nevertheless, the construction can be modified appropriately;
we will return to this question elsewhere.


\end{document}